\documentclass[11pt]{article}
\usepackage{amsfonts}
\usepackage{latexsym}

\newcommand{\bbz}{{\mathbb Z}}

\newcommand{\bbr}{{\mathbb R}}
\newcommand{\bbc}{{\mathbb C}}
\newcommand{\bbf}{{\mathbb F}}
\newcommand{\qed}{\hfill$\Box$}

\newcommand{\then}{\Longrightarrow}
\newcommand{\der}{{\Delta\,}}
\newcommand{\G}{{X}}
\newcommand{\C}{{\cal K}}
\newcommand{\T}{{\rho}}
\newcommand{\sn}{{\it span}}
\newcommand{\tr}{{\it Tr}}

\newcommand{\cont}{{\it cont}}
\newcommand{\inv}{^{-1}}
\newcommand{\inver}{{\it inv}}
\newcommand{\D}{{\it Des}}

\newcommand{\al}{\alpha}
\newcommand{\la}{\lambda}

\newcommand{\hk}{{h}}
\newcommand{\ua}{{\dot{a}}}
\newcommand{\ub}{{\dot{b}}}
\newcommand{\da}{{\ddot{a}}}
\newcommand{\db}{{\ddot{b}}}

\newtheorem{thm}{Theorem}[section]
\newtheorem{pro}[thm]{Proposition}
\newtheorem{lem}[thm]{Lemma}

\newtheorem{cor}[thm]{Corollary}

\newtheorem{exa}[thm]{Example}
\newtheorem{obs}[thm]{Observation}
\newtheorem{rem}[thm]{Remark}
\newtheorem{df}[thm]{Definition}

\begin{document}

\title{A Construction of\\ Coxeter Group Representations (II)}
\bibliographystyle{acm}
\author{Ron M. Adin%
\thanks{Department of Mathematics, Bar-Ilan University,
Ramat-Gan 52900, Israel. Email: {\tt radin@math.biu.ac.il}}\ $^\S$
\and Francesco Brenti%
\thanks{Dipartimento di Matematica,
Universit\'{a} di Roma ``Tor Vergata'',
Via della Ricerca Scientifica,
00133 Roma, Italy. Email: {\tt brenti@mat.uniroma2.it}}\ $^\S$
\and Yuval Roichman%
\thanks{Department of Mathematics, Bar-Ilan University,
Ramat-Gan 52900, Israel. Email: {\tt yuvalr@math.biu.ac.il}}
\thanks{Research of all authors was supported in part by
the Israel Science Foundation, founded by the Israel Academy of Sciences and Humanities
and
by the EC's IHRP Programme, within the Research Training Network ``Algebraic Combinatorics
in Europe'', grant HPRN-CT-2001-00272.}}
\date{submitted January 6, 2006; revised July 5, 2006\\ \medskip
Dedicated to Gordon James on the occasion of his 60th birthday}


\maketitle

\begin{abstract}
An axiomatic approach to the representation theory of Coxeter
groups and their Hecke algebras was presented in~\cite{U-I}.
Combinatorial aspects of this construction are studied in this paper.
In particular,  the symmetric group case is investigated in detail.
The resulting representations are completely classified and
include the irreducible ones.
\end{abstract}

\section{Introduction}\label{s.intro}

\subsection{Outline}
An axiomatic construction of Coxeter group representations  was
presented in~\cite{U-I}.
 This was carried out by a natural assumption on the
representation matrices, avoiding a priori use of  external
concepts (such as Young tableaux).

\smallskip

Let $(W,S)$ be a Coxeter system, and let $\C$ be a finite subset of $W$.
Let $\bbf$ be a suitable field of characteristic zero (e.g., the
field $\bbc(q)$ in the case of the Iwahori-Hecke algebra), and let
$\T$ be a representation of (the Iwahori-Hecke algebra of) $W$ on
the vector space $V_{\C}:=\sn_{\bbf}\{C_w\,|\,w\in\C\}$, with
basis vectors indexed by elements of $\C$. We want to study the
sets $\C$ and representations $\T$ which satisfy the following
axiom:

\begin{itemize}
\item[$(A)$] {\it For any generator $s\in S$ and any element $w\in
\C$ there exist scalars $a_s(w),b_s(w)\in\bbf$ such that
$$
\T_s(C_w) = a_{s}(w)C_w + b_{s}(w)C_{ws}.
$$
If $w\in\C$ but $ws\not\in\C$ we assume $b_s(w)=0$.%
}
\end{itemize}

A pair $(\T,\C)$ satisfying Axiom $(A)$ is called an {\em abstract
Young (AY) pair;} $\T$ is an {\em AY representation,} and $\C$ is
an {\em AY cell.} If $\C\ne \emptyset$ and has no proper subset
$\emptyset \subset \C' \subset \C$ such that $V_{\C'}$ is
$\T$-invariant, then $(\T,\C)$ is called a {\em minimal AY pair.}
(This is much weaker than assuming $\T$ to be irreducible.)



\medskip

In~\cite{U-I} it was shown that  an AY representation of a simply
laced Coxeter group is determined by a linear functional on the
root space. In this paper it is shown that, furthermore, the
values of the linear functional on the ``boundary" of the AY cell
determine the representation (see Theorem~\ref{t.7} below).
In Section~\ref{s.cells_Sn} this result is used
to characterize AY cells in the symmetric group.
%
This characterization is then applied to show that every
irreducible representation of
$S_n$ may be realized as a minimal abstract Young representation
(see Theorem~\ref{t.irr} below). AY representations of Weyl groups
of type $B$ are not determined by a linear functional. However, it
is shown that a similar result holds for these groups
(Theorem~\ref{B.irr} below). Finally, we characterize the elements
$\pi\in S_n$ for which the interval $[id,\pi]$ forms a minimal AY
cell, carrying an irreducible representation (see
Theorem~\ref{t.s5} below).





\subsection{Main Results}\label{s.main}




In Section~\ref{s.bound} it is shown that the action of the group
$W$ on the boundary of a cell determines the representation up to
isomorphism.

\begin{thm}\label{t.7m}{\rm (see Theorem~\ref{t.7})}
Let $(\T,\C)$ be a minimal AY pair for a simply laced Coxeter group $W$,
where $\C$ is finite.
Then the behavior of $\T$ at the boundary of $\C$
(i.e., the values $a_{s}(w)$ for $w\in \C, s\in S, ws\not\in\C$)
determines $\T$ up to isomorphism.
\end{thm}

The proof combines continuity arguments with the reduction of AY
representations to distinguished linear functionals, carried out
in~\cite{U-I} (see Theorems~\ref{t.main11} and~\ref{t.main12} below).


AY cells in the symmetric group are characterized
in Section~\ref{s.cells_Sn}.

\begin{thm}\label{t.main.symm}{\rm (see Theorem~\ref{t.symm})}
Let $\C\subseteq S_n$
and let $\sigma\in \C$. Then $\C$ is a minimal AY cell if and only
if there exists a standard skew Young tableau $Q$ of size $n$ such that
$$
\sigma^{-1}\C=\{\pi\in S_n|\ Q^{\pi\inv} \hbox{ is standard}\},
$$
where $Q^{\pi\inv}$ is the tableau obtained from $Q$ by replacing each
entry $i$ by $\pi\inv(i)$.
\end{thm}

The proof applies Theorem~\ref{t.7m} together with
Theorems~\ref{t.main11} and~\ref{t.main12} below.
Theorem~\ref{t.main.symm} is then used to prove the following.

\begin{thm}\label{t.main.irr} {\rm (see Corollary~\ref{t.irr2})}
The complete list of minimal AY representations of the symmetric
group $S_n$ is given (up to isomorphism) by the skew Specht
modules $S^{\lambda/\mu}$, where $\la/\mu$ is of order $n$
(and $\mu$ possibly empty).

In particular, every irreducible representation of the symmetric
group $S_n$ may be realized as a minimal abstract Young
representation.
\end{thm}

Combining this theorem with the combinatorial induction rule for
minimal AY representations (Theorem~\ref{induction} below) we
prove

\begin{thm}\label{B.main.irr} {\rm (see Theorem~\ref{B.irr})}
Every irreducible representation of the classical Weyl group $B_n$
may be realized as a minimal abstract Young representation.
\end{thm}





\begin{df}\label{d.main.top}
An element $w\in W$ is a {\em top} element if the interval $[id,w]$
is a minimal AY cell carrying an {\em irreducible} AY representation of $W$.
\end{df}

The top elements of the symmetric group $S_n$ are characterized in
Section~\ref{s.inter}.

\begin{thm}\label{t.main.s5} {\rm (see Theorem~\ref{t.s5})}
A permutation $\pi\in S_n$ is a top element if and only if $\pi$
is the column word of a row standard Young tableau (see
Definition~\ref{d.row} below).
\end{thm}

\noindent{\bf Note:} 
Having completed the first version of this paper, we were informed that
results equivalent to Theorems~\ref{t.7m} and~\ref{t.main.irr},
with entirely different proofs, appear in~\cite{KR, Ram03b}.

\section{Preliminaries}\label{s.prelim}

For the necessary background on Coxeter groups see~\cite{Hum};
on convex sets and generalized descent classes see~\cite{BW};
and on symmetric group representations see~\cite{Ja, JK, Sa}.
See also~\cite{DJ, KL, OV, Ram03}.


\subsection{Young Forms}

Let $Q$ be a standard Young tableau of skew shape. 
If $k\in \{1,\dots,n\}$ is in box $(i,j)$ of $Q$ 
then the {\em content} of $k$ in $Q$ is $ c(k):= j - i$.
For  $1\le k < n$, the $k$th {\em hook-distance} is
defined as $ \hk(k):=c(k+1)-c(k)$.
Denote by $Q^{s_k}$ the tableau obtained from $Q$ by interchanging
$k$ and $k+1$. 
The classical Young orthogonal form for $S_n$ (see, e.g.,
\cite[\S 25.4]{Ja}) is generalized naturally to skew shapes.

\begin{thm}\label{t.yof-skew}
{\bf (Young Orthogonal Form for Skew Specht Modules)}
Let $\{v_Q\,|\,Q$ standard Young tableau of shape $\la/\mu\}$
be the basis of the skew Specht module $S^{\la/\mu}$ obtained by
the Gram-Schmidt process from the polytabloid basis. Then
\begin{equation}\label{e.yofi}
\T^{\la/\mu}(s_i)(v_Q)= \frac{1}{\hk(i)}\,v_Q +
\sqrt{1-\frac{1}{\hk(i)^2}}\,v_{Q^{s_i}}.
\end{equation}
\end{thm}

\noindent{\bf Proof.}
(Due to J.\ Stembridge~\cite{Ste}; see also~\cite{Greene}.)
Matrices determined by (\ref{e.yofi})
must satisfy the Coxeter relations of $S_n$, because the same is true
when the skew tableaux are completed to full tableaux of non-skew shape.
Therefore they define a representation of $S_n$, which we denote
$Y^{\lambda/\mu}$. Upon restricting the action of $S_n$ to
$S_k\times S_{n-k}$ (where $n=|\lambda|, k=|\mu|$), $Y^\lambda$ decomposes
into the direct sum
$\bigoplus_{\{\mu\subseteq\lambda\,|\,|\mu|=k\}}Y^\mu \otimes Y^{\lambda/\mu}$.
On the other hand, Specht modules have exactly the same decomposition.
This follows, for example, from the corresponding identity on Schur
functions~\cite[(7.66)]{St2} (using the inverse Frobenius image).
Since $Y^\lambda\cong S^\lambda$, $Y^{\lambda/\mu}$ must be
isomorphic to $S^{\lambda/\mu}$.
\qed

\bigskip

$B_n$, the classical Weyl group of type $B$, is a Coxeter system
with $S=\{s_i|\ 0\le i<n\}$, $m(s_0,s_1)=4$,
$m(s_i,s_{i+1})=3$ for $1\le i<n$, and $m(s_i,s_j)=2$ otherwise.
The irreducible representations of $B_n$ are indexed by
pairs of partitions $(\lambda,\mu)$, where
$\lambda$ is a partition of some $0\le k\le n$ and
$\mu$ is a partition of $n-k$.
A basis for the irreducible representation of shape $(\lambda,\mu)$
may be indexed by all pairs $(P,Q)$ of standard Young tableaux of shapes
$\lambda$ and $\mu$, respectively, where
$P$ is a tableau on a subset of $k$ letters from $\{1,\dots,n\}$ and
$Q$ is a tableau on the complementary subset of letters.
There exists a basis such that the following Young form holds
(see, e.g., \cite{P}).

\begin{thm}\label{B-yof}{\bf (Classical Young Orthogonal Form for $B_n$)}
Denote the above basis elements by $v_{(P,Q)}$.
For $1\le i<n$ define the hook distance $\hk(i)$ as follows:
$$
\hk(i) :=
\cases{%
h_P(i), &if $i$ and $i+1$ are both in $P$;\cr
h_Q(i), &if $i$ and $i+1$ are both in $Q$;\cr
\infty, &if $i$ and $i+1$ are in different tableaux.}
$$
Then, for $1\le i<n$,
$$
\rho^{\la,\mu}(s_i)(v_{(P,Q)}) =
\frac{1}{\hk(i)}\,v_{(P,Q)} +
\sqrt{1-\frac{1}{\hk(i)^2}}\,v_{{(P,Q)}^{s_i}},
$$
where ${(P,Q)}^{s_i}$ is the pair of tableaux obtained from
$(P,Q)$ by interchanging $i$ and $i+1$, whereas
$$
\rho^{\la,\mu}(s_0)(v_{(P,Q)})=\cases{%
v_{(P,Q)}, &if $1$ is in $P$;\cr
-v_{(P,Q)}, &if $1$ is in $Q$.}
$$
\end{thm}

\subsection{Abstract Young Representations}

Recall the definition of AY cells and representations from the
introduction.

\begin{pro}\label{t.main.convex}{\rm\cite[Corollary 4.4]{U-I}}
Every minimal AY cell is convex in the Hasse diagram of the
right weak Bruhat order.
\end{pro}

\begin{df}
For a convex subset $\C\subseteq W$ define:
\begin{eqnarray*}
T_{\C}          &:=& \{wsw^{-1}\,|\,s\in S,\,w\in \C,\,ws\in \C\},\\
T_{\partial \C} &:=& \{wsw^{-1}\,|\,s\in S,\,w\in \C,\,ws\not\in \C\}.
\end{eqnarray*}
\end{df}

\begin{df}\label{d.cg} ($\C$-genericity)\\
Let $\C$ be a convex subset of $W$ containing the identity element.
A linear functional $f$ on the root space $V$ is {\rm $\C$-generic} if:
\begin{itemize}
\item[{\bf (i)}]
For all $t\in T_{\C}$,
$$
\langle f,\al_t\rangle \not\in \{0, 1, -1\}.
$$
\item[{\bf (ii)}]
For all $t\in T_{\partial \C}$,
$$
\langle f,\al_t\rangle \in \{1, -1\}.
$$
\item[{\bf (iii)}]
If $w\in \C$, $s,t\in S$, $m(s,t)=3$ and $ws, wt\not\in \C$ then
$$
\langle f,\al_{wsw^{-1}}\rangle =
\langle f,\al_{wtw^{-1}}\rangle\;(= \pm1).
$$
\end{itemize}
\end{df}


By~\cite[Observation 3.3]{U-I}, we may assume that
$id\in\C$. By~\cite[Theorem 11.1]{U-I}, under mild conditions,
Axiom $(A)$ is equivalent to the following:

\begin{itemize}
\item[$(B)$] {\it For any reflection $t$ there exist scalars
$\ua_t, \ub_t, \da_t, \db_t \in\bbf$ such that, for all $s\in S$
and $w\in \C$:
$$
\T_s(C_w) =\cases{%
\ua_{wsw^{-1}}C_w+\ub_{wsw^{-1}}C_{ws}, &if $\ell(w)<\ell(ws)$;\cr
\da_{wsw^{-1}}C_w+\db_{wsw^{-1}}C_{ws}, &if
$\ell(w)>\ell(ws)$.\cr}
$$
}
\end{itemize}

\begin{thm}\label{t.main11}{\rm\cite[Theorem 7.4]{U-I}}
Let $(W,S)$ be an irreducible simply laced Coxeter system,
and let $\C$ be a convex subset of $W$ containing the identity element.
If $f\in V^*$ is $\C$-generic then
$$
\ua_t := {1\over \langle f,\al_t\rangle}
\qquad(\forall t\in T_{\C} \cup T_{\partial\C}),
$$
together with $\da_t$, $\ub_t$ and $\db_t$ satisfying
\begin{eqnarray*}
\ua_t + \da_t &=& 0\\
\ub_t \cdot \db_t &=& (1 - \ua_t)(1 - \da_t)
\end{eqnarray*}
define a representation $\T$ such that $(\T,\C)$ is a minimal AY pair
satisfying Axiom $(B)$.
\end{thm}

The following theorem is complementary.

\begin{thm}\label{t.main12}{\rm\cite[Theorem 7.5]{U-I}}
Let $(W,S)$ be an irreducible simply laced Coxeter system
and let $\C$ be a subset of $W$ containing the identity element.
If $(\T,\C)$ is a minimal AY pair satisfying Axiom $(B)$
and $\ua_t\ne 0$ ($\forall t\in T_{\C}$)
then
there exists a $\C$-generic $f\in V^*$ such that
$$
\ua_t = \frac{1}{\langle f,\al_t\rangle}\qquad
(\forall\ t\in T_{\C} \cup T_{\partial\C}).
$$
\end{thm}




The following combinatorial rule for induction of AY
representations is analogous to the one for Kazhdan-Lusztig
representations~\cite{BV2,Geck}.


\begin{thm}\label{induction}{\rm \cite[Theorem 9.3]{U-I}}
Let $(W,S)$ be a finite Coxeter system, $P=\langle J\rangle$
$(J\subseteq S)$ a parabolic subgroup,
and $W^J$ the set of all representatives of minimal length of the
right cosets of $P$ in $W$. Let $(\psi,\cal D)$ be a minimal AY
pair for $P$. Then
\begin{itemize}
\item[1.]
${\cal D} W^J$ is a minimal AY cell for $W$.
\item[2.]
The induced representation $\psi\uparrow ^W_P$ is isomorphic to
an AY representation on $V_{{\cal D} W^J}$.
\end{itemize}
\end{thm}

\begin{rem}\label{induction1}
By~\cite[Lemma 9.7]{U-I}, for every $s\in S$ and $r\in W^J$
either $rs\in W^J$, or $rs\not\in W^J$ and $rs=pr$ with $p\in J$.
Then, by the proof of~\cite[Theorem 9.3]{U-I},
the representation matrices of the generators in the
resulting induced representation are as follows:
for $s\in S$, $m\in {\cal D}$, $r\in W^J$,
$$
\rho_s(C_{mr}) = \cases{%
C_{mrs}, &if $rs\in W^J$;\cr a_p(m)C_{mr} + b_p(m)C_{mrs},
&otherwise $(rs = pr,\,p\in J)$,}
$$
where the coefficients $a_p$ and $b_p$ are given by the AY
representation $\psi$; namely, $ \psi_p(C_m) = a_p(m) C_m + b_p(m)
C_{mp}.$
\end{rem}

\section{Boundary Conditions}\label{s.bound}

In this section it is shown that the action of the group $W$ on
the boundary of a minimal AY cell determines the representation up to isomorphism.







For a subset of reflections $A$ let the (left) {\em $A$-descent
set} of an element $w\in W$ is defined by
$$
\D_A(w):=\{t\in A\,|\,\ell(tw)<\ell(w)\}.
$$

\begin{df}\label{d.CF}
Let $w\in W$, and let $f\in V^*$ be an arbitrary linear functional on
the root space $V$ of $W$.
\begin{itemize}
\item[(1)] Define
$$
A_f := \{t\in T\,|\,\langle f,\al_t\rangle \in\{1, -1\}\},
$$
and
$$
\C^f_{w}:=\{v\in W|\D_A(v)=\D_A(w)\}.
$$
\item[(2)] If $f$ is $\C^f_w$-generic (as in Definition~\ref{d.cg})
then the corresponding AY representation of $W$ (as in Theorem~\ref{t.main11}),
with the
symmetric normalization $\db_t =\ub_t\,(\forall t\in T_{\C^f_w})$,
will be denoted $\T^f_w$ (or just $\T^f$ in case there is no ambiguity).
\end{itemize}
\end{df}


\begin{rem}\label{r.indep}
By~\cite[Theorem 11.1]{U-I}, the representation $\T^f_w$
is independent of the normalization (up to isomorphism).
\end{rem}


\begin{df}\label{d.flat}
Let $W$ be a Coxeter group, and let $V$ be its root space.
A {\em basic (affine) hyperplane} in $V^*$ has the form
$$
H_{t,\varepsilon}:=\{f\in V^* \,|\, \langle
f,\al_t\rangle=\varepsilon\},
$$
for some $t\in T$ and $\varepsilon \in \{1,-1\}$.

A {\em basic flat} in $V^*$ is an intersection of basic hyperplanes.
It is {\em proper} if different from $\emptyset$ and $V^*$.

For a basic flat $L$, let
$$
A = A_L := \{\,t\in T \,|\, L\subseteq H_{t,\varepsilon}
\hbox{\rm\ for some\ } \varepsilon \in \{1,-1\}\,\}.
$$
Then $\{W_A^D \,|\, D\subseteq A\}$, where $W_A^D:=\{ w\in
W|\D_A(w)=D\}$,
is a partition of $W$ into convex
subsets, called the {\em $L$-partition} of $W$.
\end{df}


Note that, for the two ``improper'' flats:
\begin{eqnarray*}
L = \emptyset &\then& A_L = T\\
L = V^*       &\then& A_L = \emptyset
\end{eqnarray*}

\begin{thm}\label{t.5}
Let $W$ be a simply laced Coxeter group. Let $L$ be a basic 
flat in $V^*$, and fix a nonempty finite convex set $\C$ in the
$L$-partition of $W$. Then, for any two elements $v,v'\in\C$ and
any two $\C$-generic vectors $f,f'\in L$,
$\C^f_v = \C^{f'}_{v'} = \C$ and
the representations $\T^f_v$ and $\T^{f'}_{v'}$ are isomorphic.
\end{thm}

\noindent{\bf Proof.}
First of all,
$$
f\in L \iff \langle f, \al_t \rangle = \varepsilon_t\quad(\forall t\in A_L)
       \iff A_L \subseteq A_f
$$
and therefore, for any $\C$ in the $L$-partition of $W$ and any $v\in\C$,
$$
\C^f_v \subseteq \C.
$$
If $f$ is also $\C$-generic then $\langle f, \al_t \rangle \ne \pm1$
for all $t\in T_{\C}$, so that $\C^f_v = \C$.

Now choose $f_0\in L$, and let $\{f_1,\dots,f_k\}$ be a basis for
the linear subspace $L-f_0$ of $V^*$. Each $f\in L$ has a unique expression as
$$
f=f_0+r_1 f_1+\dots+r_k f_k,
$$
where $r_1,\ldots,r_k\in \bbr$. For any $t\in T_{\C} \cup T_{\partial\C}$,
$\langle f,\al_t\rangle$ is a linear combination of $1,r_1,\ldots,r_k$,
and is nonzero if $f$ is $\C$-generic.
For $v\in\C$, use the represetation $\T^f_v$ with the row-stochastic 
normalization $\ua_t + \ub_t = \da_t + \db_t = 1\,(\forall t\in T_{\C^f_v})$;
see Remark~\ref{r.indep}.

Thus, for any $v\in\C$ and $\C$-generic $f\in L$, each entry of the matrix
$\T^f_v(s)$ ($\forall s\in S$) is a rational function of $r_1,\dots,r_k$;
and the same therefore holds for each entry of $\T^f_v(w)$ ($\forall w\in W$)
and for the character values $\tr(\T^f_v(w))$. Note that these
rational functions (unlike the actual values of $r_1,\ldots,r_k$)
do not depend on the choice of $v$ and $f$, even though
the set $L^{\rm gen}$ of all $\C$-generic $f\in L$ may be disconnected
(see example below). By discreteness of the character values and
continuity of the rational function, each character value is
constant in each connected component of $L^{\rm gen}$,
and at the same time represented by one rational function
throughout $L^{\rm gen}$. It is therefore the same constant for all
$f\in L^{\rm gen}$ (and $v\in\C$), as claimed.
\qed

\begin{exa}
Take $W=S_3=\langle s_1,s_2\rangle$ (type $A_2$) and the basic
flat $L = \{f\in V^*\ |\ \langle f,\al_{s_1s_2s_1}\rangle=-1\}$.
Then $A = \{s_1s_2s_1\}$, and we may choose $\C =\{id,s_1,s_2\}$.
In that case, $T_{\C} = \{s_1,s_2\}$ and
$T_{\partial\C} = \{s_1s_2s_1\} = A$.
$L$ is an affine line in $V^* \cong \bbr^2$, and the $\C$-generic points
in $L$ form five disjoint open intervals (three of them bounded). For any
$\C$-generic vector $f\in L$ and any $v\in\C$, $\T^f_v$ is the 3-dimensional
representation isomorphic to the direct sum of the sign representation and
the unique irreducible 2-dimensional representation of $S_3$.
\end{exa}

An important special case is $L = V^*$ ($\C = W$).

\begin{thm}\label{t.1}
Let $W$ be a finite simply laced Coxeter group, and let $f\in V^*$
be $W$-generic (i.e.,
$\langle f,\al_t\rangle \not\in\{0,1,-1\}, \forall t\in T$).
Then, for any $v\in W$, the representation $\rho^f_v$ on $V_W$
is isomorphic to the regular representation of $W$.
\end{thm}

\noindent{\bf Proof.}
Fix $v\in W$ (and ignore it in the notation).
For all but finitely many values of $\mu\in \bbr$, the
linear functional $\mu f\in V^*$ is also $W$-generic.
The representations $\rho^{\mu f}$ and $\rho^f$ are isomorphic,
by Theorem~\ref{t.5}. On the other hand, if $|\mu| \to \infty$ then
$$
a_s^{(\mu f)}(w) = {\pm1\over \langle \mu f,\al_{wsw^{-1}}\rangle}\to 0
\qquad (\forall s\in S, w\in W)
$$
and consequently $b_s^{(\mu f)}(w)\to 1$. The representation
matrices of $\rho^{(\mu f)}(s)$ $(\forall s\in S)$, and thus
also those of $\rho^{(\mu f)}(w)$ $(\forall w\in W)$, tend to those of
the regular representation. The character of $\rho^f$ is thus the
character of the regular representation.
\qed

\medskip

Theorem~\ref{t.5} may be reformulated as follows.

\begin{thm}\label{t.7}
Let $(\T,\C)$ be a minimal AY pair for a simply laced Coxeter group $W$,
where $\C$ is finite.
Then the behavior of $\T$ at the boundary of $\C$
(i.e., the values $a_{s}(w)$ for $w\in \C, s\in S, ws\not\in\C$)
determines $\T$ up to isomorphism.
\end{thm}

%

\section{Minimal Cells in $S_n$}\label{s.cells_Sn}

In this section we show that integer-valued $\C$-generic vectors for $W=S_n$
lead to standard Young tableaux (of skew shape).
Theorem~\ref{t.5} is then applied to give a
complete characterization of minimal AY cells in $S_n$. Finally,
it is shown that all irreducible representations of $S_n$ are
minimal AY.


\subsection{Identity Cells and Skew Shapes}\label{s.IdCT}

In this subsection we study minimal AY cells $\C\subseteq S_n$. 
By~\cite[Observation 3.3]{U-I}, every minimal AY cell is
a translate of a minimal AY cell containing the identity element;
thus we may assume that $id\in\C$.

\medskip

For a vector $v=(v_1,\dots,v_n)\in \bbr^n$ denote
$$
\der v:=(v_2-v_1,\dots,v_n-v_{n-1})\in \bbr^{n-1}.
$$
For a (skew) standard Young tableau $Q$ denote $c(k):=j-i$, where $k$ is
the entry in row $i$ and column $j$ of $Q$.
Call $cont(Q):=(c(1),\dots,c(n))$ the {\em content vector} of $Q$, and
call $\der cont(Q)$ the {\em derived content vector} of $Q$.

Note that for $W=S_n$ we may identify the root space $V$ with a subspace
(hyperplane) of $\bbr^n$:
$$
V \cong \{(v_1,\ldots,v_n)\in\bbr^n\,|\,v_1+\ldots+v_n=0\}.
$$

\noindent
The positive root $\al_{ij}\in V$ corresponding to
the transposition $(i,j)\in S_n$ may be identified with
the vector $\varepsilon_i - \varepsilon_j$ ($1\le i < j\le n$),
where $\{\varepsilon_1,\ldots,\varepsilon_n\}$ is
the standard basis of $\bbr^n$.
The dual space $V^*$ is then a {\em quotient} of $\bbr^n$:
$$
V^* \cong \bbr^n/\,\bbr e,
$$
where $e := (1,\ldots,1)\in \bbr^n$.
We shall abuse notation and represent a linear functional $f\in V^*$ by
any one of its representatives $f=(f_1,\ldots,f_n)\in \bbr^n$; the
natural pairing $\langle\cdot,\cdot\rangle : V^* \times V \to \bbr$ is
then given by $\langle f,\varepsilon_i-\varepsilon_j\rangle = f_i-f_j$.


\medskip

Recall the notations $\C^f_w$ and $\T^f_w$ from Definition~\ref{d.CF}.


\begin{thm}\label{t.tabl2}
Let  $f\in \bbr^n$ have integer coordinates. Then:
$(\T^f_{id},\C^f_{id})$ is a minimal AY pair for $W=S_n$ if and only if
there exists a standard skew Young tableau $Q$ of size $n$ such that
$$
\der f=\der cont(Q).
$$
\end{thm}


The proof of Theorem~\ref{t.tabl2} relies on the following lemmas.

%

\begin{lem}\label{t.monotone}
Let $f\in \bbr^n$ and $1\le i < j\le n$.
If either $\langle f,\al_{ij}\rangle = \pm 1$,
or $f$ is $\C^f_{id}$-generic and $\langle f,\al_{ij}\rangle = 0$,
then $w^{-1}(i)<w^{-1}(j)$ for all $w\in \C^f_{id}$.
\end{lem}

\noindent{\bf Proof.}
The claim clearly holds for $w=id$. It thus suffices to show that
if $w,ws\in \C^f_{id}$ ($s\in S$) then $w^{-1}(j)-w^{-1}(i)$ and
$(ws)^{-1}(j)-(ws)^{-1}(i)$ have the same sign.

Since $s$ is an adjacent transposition, say $s=(t,t+1)$ ($1\le t\le n-1$),
the two signs differ if and only if $\{w^{-1}(i),w^{-1}(j)\}=\{t,t+1\}$.
This implies that $(i,j)=wsw^{-1}\in T_{\C^f_{id}}$.
Thus $\langle f,\al_{ij}\rangle \ne\pm 1$ and, if $f$ is $\C^f_{id}$-generic,
also $\langle f,\al_{ij}\rangle \ne 0$. This contradicts the assumption.
\qed

\begin{lem}\label{t.rep-cell2}
Let $f\in \bbr^n$ be an arbitrary vector.
Then $f$ is $\C^f_{id}$-generic if and only if, for all $1\le i < j\le n$:
\begin{equation}\label{e.R1}
\langle f, \al_{ij}\rangle=0\ \then\
\exists\ r_1,r_2\in[i+1,j-1]\ \hbox{\rm s.t.}\
\langle f, \al_{ir_1}\rangle=-\langle f, \al_{ir_2}\rangle=1.
\end{equation}
\end{lem}

\noindent{\bf Proof.}
Let $\C := \C^f_{id}$.

\noindent{\bf A - (necessity).}
%
Note that, since $f$ is $\C$-generic,
\begin{equation}\label{e.9}
\langle f, \alpha_{ij}\rangle=0 \then (i,j)\not\in T_\C\cup T_{\partial \C}.
\end{equation}
Consider the set
$$
Z:=\{(i,j)\ |\ 1\le i<j\le n, \langle f,\al_{ij}\rangle =0\}.
$$
We shall prove that
$$
(i,j)\in Z \then \exists\,r_1,r_2\in[i+1,j-1] \hbox{ such that }
\langle f,\al_{ir_1}\rangle=\langle f,\al_{r_2 j}\rangle=1.
$$
The proof will proceed by induction on $j-i$,
the height of the root $\al_{ij}$.

\smallskip

Assume first that $j-i = 1$. Then $(i,j)=(i,i+1)\in S$.
Since $id\in\C$, $S \subseteq T_{\C} \cup T_{\partial\C}$.
This contradicts (\ref{e.9}) above.



For the induction step, assume
that $(i,j)\in Z$ with $j-i>1$ and that the claim is true for all
reflections in $Z$ with smaller heights.
Choose $w\in \C$ such that $d_w:=|w^{-1}(j)-w^{-1}(i)|$ is minimal.
Note that, by Lemma~\ref{t.monotone}, actually $d_w=w^{-1}(j)-w^{-1}(i)>0$.

If $d_w=1$ then there exists $1\le t\le n-1$ such that $w(t)=i$ and $w(t+1)=j$,
so that $(i,j) = w(t,t+1)w^{-1} \in T_{\C} \cup T_{\partial\C}$,
which is a contradiction to (\ref{e.9}).

Thus $d_w\ge 2$.

Define $1\le r_1,r_2\le n$ by $w^{-1}(r_1)=w^{-1}(i)+1$ and $w^{-1}(r_2)=
w^{-1}(j)-1$. By minimality of $d_w$, $(i,r_1)\in T_{\partial\C}$ so that
$\langle f,\al_{i r_1}\rangle = \pm1$; similarly
$\langle f,\al_{r_2 j}\rangle = \pm1$.
Now $\langle f, \al_{i j}\rangle = 0$ and $\langle f, \al_{i r_1}\rangle = \pm1$
imply
\begin{eqnarray*}
\langle f, \al_{r_1 j}\rangle &=&
\langle f, \al_{r_1 i}\rangle + \langle f, \al_{i j}\rangle = \pm1\quad
\hbox{\rm (if $r_1 < i$)};\cr
\langle f, \al_{r_1 j}\rangle &=&
\langle f, \al_{i j}\rangle - \langle f, \al_{i r_1}\rangle = \pm1\quad
\hbox{\rm (if $i < r_1 < j$)};\cr
\langle f, \al_{j r_1}\rangle &=&
\langle f, \al_{i r_1}\rangle - \langle f, \al_{i j}\rangle = \pm1\quad
\hbox{\rm (if $j < r_1$)}.
\end{eqnarray*}
Since $w^{-1}(i) < w^{-1}(r_1) < w^{-1}(j)$ we conclude, by Lemma~\ref{t.monotone},
that $i < r_1 < j$. Similarly  $i < r_2 < j$.

If $\langle f,\al_{i r_1}\rangle = -\langle f,\al_{i r_2} \rangle = 1$ or
$\langle f,\al_{i r_1}\rangle = -\langle f,\al_{i r_2} \rangle = -1$
we are done.
We can thus assume, with no loss of generality, that
$\langle f,\al_{i r_1}\rangle = \langle f,\al_{i r_2} \rangle = \varepsilon = \pm1$
and $r_1\le r_2$.

If $r_1=r_2$ then $w^{-1}(j)-w^{-1}(i)=2$. Denote $t := w^{-1}(i)$.
Then
$\langle f,\al_{w (t,t+1)w^{-1}}\rangle = \langle f,\al_{ir_1}\rangle = \varepsilon$
and $\langle f,\al_{w (t+1,t+2)w^{-1}}\rangle = \langle f,\al_{r_1 j}\rangle = -\varepsilon$.
This contradicts condition (iii) of $\C$-genericity (Definition~\ref{d.cg}).
Therefore $r_1<r_2$.
Thus
$$
\langle f,\al_{r_1 r_2}\rangle =
\langle f,\al_{i r_2}\rangle - \langle f,\al_{i r_1}\rangle = 0.
$$
Since $r_2-r_1<j-i$, by the induction hypothesis there exists
$r_1<r_3<r_2$ such that $\langle f,\al_{r_1 r_3}\rangle=-\varepsilon$.
Thus $\langle f,\al_{i r_3} \rangle=0$
(and $\langle f,\al_{r_3 j} \rangle=0$).
Again, by the induction hypothesis, there exist
$i<r_4,r_5<r_3$ such that
$\langle f,\al_{i r_4} \rangle=
\langle f,\al_{r_5 r_3} \rangle=1$.
Noting that
$\langle f,\al_{r_5 r_3} \rangle=
\langle f,\al_{r_5 j} \rangle$ completes the proof that
condition~(\ref{e.R1}) is necessary.

\medskip

\noindent{\bf B - (sufficiency).}

Assume now that $f\in \bbr^n$ satisfies condition~(\ref{e.R1}).
Condition (ii) of Definition~\ref{d.cg} holds by the definition of $\C^f_{id}$.
Assume that $\langle f,\al_{i j} \rangle=0$.
By condition~({\ref{e.R1}) and Lemma~\ref{t.monotone},
there exist
$i<r_1<r_2<j$ such that $w^{-1}(i)<w^{-1}(r_1)<w^{-1}(r_2)<w^{-1}(j)$
for all $w\in \C$.
Thus $w^{-1}(j)-w^{-1}(i)>2$, and this implies conditions (i) and (iii)
of $\C$-genericity as follows :

For condition (i), if $w,ws\in \C$, $s=(t,t+1)\in S$
and $\langle f,\al_{wsw^{-1}} \rangle=0$ then, denoting $i:=w(t)$ and $j:=w(t+1)$,
we get $i<j$ and $\langle f,\al_{i j} \rangle=0$, so that
$w^{-1}(j)-w^{-1}(i)=1$ contradicting our conclusion above.

For condition (iii), if $w\in \C$, $ws,wt\not\in\C$,
$s=(k,k+1)$ and $t=(k+1,k+2)$
then denote $i:=w(k)$, $r:=w(k+1)$, and $j:=w(k+2)$. Then $i<r<j$ and
$\langle f,\al_{wsw^{-1}} \rangle=\pm 1$, $\langle
f,\al_{wtw^{-1}} \rangle=\pm 1$. If $\langle f,\al_{wsw^{-1}}
\rangle \ne \langle f,\al_{wtw^{-1}} \rangle$ then $\langle
f,\al_{wsw^{-1}} \rangle + \langle f,\al_{wtw^{-1}} \rangle =0$,
that is $\langle f,\al_{ir} \rangle + \langle f,\al_{rj} \rangle
=0$ or equivalently $\langle f,\al_{ij} \rangle=0$, and
$w^{-1}(j)-w^{-1}(i) = 2$ contradicts our conclusion above. \qed

\medskip

\begin{lem}\label{t.rep-cell4}
A vector  $c=(c_1,\ldots,c_n)\in \bbz^n$ is a content
vector for some skew standard Young tableau if and only if for all $1\le i<j\le n$
\begin{equation}\label{e.cont}
c_i = c_j \then \exists\,r_1, r_2\in[i+1,j-1] \hbox{ such that }
c_{r_1} = c_i+1 \hbox{ and } c_{r_2} = c_i-1.
\end{equation}
\end{lem}

\noindent{\bf Proof.}
It is clear that if $(c_1,\ldots,c_n)$ is the content vector of
a skew standard Young tableau then it satisfies condition~(\ref{e.cont}).

Conversely, let $(c_1,\ldots,c_n)\in \bbz^n$ be such that (\ref{e.cont}) holds.
We will show, by induction on $n$, that there exists
a skew standard Young tableau $Q$ such that $\cont(Q) = (c_1,\ldots,c_n)$.
The existence of $Q$ is clear for $n\le 2$.
By the induction hypothesis, there exists a skew standard Young tableau
$Q'$ such that $\cont(Q') = (c_1,\ldots,c_{n-1})$.
Let $C := \{c_\ell\,|\,\ell\in[n-1]\}$.

If $c_n\in C$, let
$$
k := \max \{\ell\in [n-1]: \; c_{\ell}=c_{n}\}.
$$
By our hypothesis there exist $r_1,r_2 \in [k+1,n-1]$ such that
$c_{r_1}=c_k + 1$ and $c_{r_2}=c_k - 1$.
If $k$ is in box $(i,j)$ of $Q'$ then
box $(i+1,j+1)$ must be empty (since $k$ is maximal). Therefore
$r_1$ must be in box $(i,j+1)$ and $r_2$ must be in box $(i+1,j)$.
Placing $n$ in box $(i+1,j+1)$ yields a skew standard Young tableau $Q$
such that $\cont(Q) = (c_1,\ldots,c_n)$, as desired.

If $c_n\not\in C$ then $Q'$ is the disjoint union of two (possibly empty)
tableaux, $Q'_+$ and $Q'_-$, consisting of the boxes of $Q'$ with contents
strictly larger (respectively, smaller) than $c_n$.
Let $(i_+,j_+)$ be the (unique) box with the smallest (closest to $c_n$)
content in $Q'_+$, and define similarly $(i_-,j_-)$ for $Q'_-$.
All of $Q'_+$ is (weakly) northeast of $(i_+,j_+)$,
all of $Q'_-$ is (weakly) southwest of $(i_-,j_-)$, and
$(i_+,j_+)$ is (strictly) northeast of $(i_-,j_-)$.
If the difference in contents between $(i_+,j_+)$ and $(i_-,j_-)$ is $2$
(the smallest possible) then these boxes have a common corner, and we can
place $n$ in box $(i_-,j_- +1) = (i_+ +1, j_+)$ to form $Q$.
If the difference is larger than we can shift all the boxes of $Q'_+$
diagonally (preserving their contents) until $i_+ = i_- -1$ (and thus
$j_+ > j_- +1$). Now we can place $n$ in box $(i_-,j_- +1)$ to form $Q$.
The discussion is even simpler if either one of $Q'_+$ and $Q'_-$ is empty.
\qed

\bigskip

\noindent{\bf Proof of Theorem~\ref{t.tabl2}.} Combine
Lemma~\ref{t.rep-cell2} with Lemma~\ref{t.rep-cell4}. \qed

\subsection{Cell Elements and Standard Tableaux}\label{s.CET}

By Theorem~\ref{t.tabl2}, a minimal AY cell (containing the identity) in $S_n$
is defined by a linear functional represented by a vector $f\in \bbz^n$
such that $\der f=\der cont(Q)$ for some standard skew Young tableau $Q$.
We will show that there is a bijection between the elements of $\C^f_{id}$
and the standard Young tableaux of the same shape as $Q$.



\begin{thm}\label{t.cell}
Let $Q$ be a standard skew Young tableau, and let $f\in \bbz^n$
be any vector satisfying $\der f = \der cont(Q)$.
Then, for any $\pi\in S_n$,
$$
\pi\in \C^f_{id} \iff \hbox{\ the tableau }Q^{\pi\inv}\hbox{\ is standard},
$$
where $Q^{\pi\inv}$ is the tableau obtained from $Q$
by replacing each entry $i$ by $\pi\inv(i)$ ($1\le i\le n$).
\end{thm}

\begin{cor}\label{t.card}
The size of $\C^f_{id}$ is equal to the number of standard Young tableaux
of the same shape as $Q$.
\end{cor}


In order to prove Theorem~\ref{t.cell}, we first make the following observation.

\begin{obs}\label{t.adjacent}
For a standard skew Young tableau $Q$ and any $1\le i< n$, exactly one of
the following 3 cases holds:
\begin{itemize}
\item[(1)]
$i+1$ is adjacent to $i$ in the same row of $Q$, and then
$$
[\der cont(Q)]_i = 1.
$$
\item[(2)]
$i+1$ is adjacent to $i$ in the same column of $Q$, and then
$$
[\der cont(Q)]_i = -1.
$$
\item[(3)]
$i+1$ and $i$ are not in the same row or column of $Q$, and then
$$
|[\der cont(Q)]_i| \ge 2.
$$
\end{itemize}
Note that $i+1$ and $i$ cannot be in the same diagonal of $Q$:
$[\der cont(Q)]_i \ne 0$.
\end{obs}

\begin{lem}\label{t.si}
Assume that $\pi\in \C^f_{id}$ and $Q^{\pi\inv}$ is standard. Then, for any
$1\le i< n$:
$$
\pi s_i\in \C^f_{id} \iff Q^{(\pi s_i)\inv} \hbox{\ is standard}.
$$
\end{lem}

\noindent{\bf Proof.}
Consider $\pi\in \C^f_{id}$ and $1\le i< n$. Then:
\begin{eqnarray}\label{e.78}
\frac{1}{\ua_{\pi s_i \pi^{-1}}} &=&
\langle f,\alpha_{\pi s_i \pi^{-1}} \rangle =
\langle f,\alpha_{(\pi(i),\pi(i+1))} \rangle =\cr
&=& \pm\,(f_{\pi(i+1)} - f_{\pi(i)}) =
\pm\,[\der cont(Q^{\pi\inv})]_i,
\end{eqnarray}
where ``$\pm$'' is the sign of $\pi(i+1) - \pi(i)$.
Thus
$$
\ua_{\pi s_i \pi^{-1}} \ne \pm 1 \iff
[\der cont(Q^{\pi\inv})]_i \ne \pm 1.
$$
On the other hand, since $\pi\in\C^f_{id}$,
$$
\ua_{\pi s_i\pi^{-1}}\ne \pm 1 \iff \pi s_i\in \C^f_{id}.
$$
By Observation~\ref{t.adjacent}, this means that
$\pi s_i\in \C^f_{id}$ if and only if $i$ and $i+1$ are not in
the same row or column of $Q^{\pi\inv}$.
Thus, for $\pi\in \C^f_{id}$ with $Q^{\pi\inv}$ standard:
$$
\pi s_i\in \C^f_{id} \iff
(Q^{\pi\inv})^{s_i} = Q^{s_i \pi\inv} = Q^{(\pi s_i)\inv}
\hbox{\ is standard}.
$$
\qed

\medskip

\noindent{\bf Proof of Theorem~\ref{t.cell}.}
By Lemma~\ref{t.si}, it suffices to show that any $\pi\in\C^f_{id}$
may be reduced to the identity permutation by
a sequence of multiplications (on the right) by adjacent transpositions
$s_i\in S$ such that all the intermediate permutations are also in $\C^f_{id}$;
and that a similar property holds for any $\pi\in S_n$ such that
$Q^{\pi\inv}$ is standard. In other words, we need to show that
$\C^f_{id}$ and $\{\pi\in S_n\,|\,Q^{\pi\inv}\hbox{\ is standard}\}$ are
{\em connected} subsets in the right Cayley graph of $S_n$ with respect to
the Coxeter generators.

For $\C^f_{id}$ this follows from the {\em convexity} of minimal
AY cells 
(Proposition~\ref{t.main.convex}).

For $\{\pi\in S_n\,|\,Q^{\pi\inv}\hbox{\ is standard}\}$ we give the
outline of an argument.
An {\em inversion} in a standard skew Young tableau $Q$ is a pair $(i,j)$
such that $1\le i < j \le n$ and $i$ appears in $Q$ strictly south of $j$.
The {\em inversion number} $\inver(Q)$ is the number of inversions in $Q$
(see~\cite{St_sign}).
If $i$ appears in $Q$ strictly south of $i+1$ Then $Q^{s_i}$ is also
a standard tableau, with $\inver(Q^{s_i}) = \inver(Q)-1$. Thus
every standard tableau $Q$ leads, by a sequence of applications of
generators $s_i\in S$, to the unique standard tableau of the same skew shape
for which $i$ is always weakly north of $i+1$ ($1\le i < n$), i.e.,
the corresponding {\em row tableau} (see Definition~\ref{d.row} below).
Thus any two standard skew tableaux of the same shape are connected by
such a sequence, and this is the connectivity result that we need.
\qed

\bigskip

%

%

In contrast to Kazhdan-Lusztig theory,
where the bijection between cell elements and tableaux
is given by the RSK algorithm,
the above bijection between elements of the cell $\C^f_{id}$ and tableaux
is extremely simple.

%
%
%



\bigskip

A complete characterization of minimal AY cells in $S_n$ now follows.

\begin{thm}\label{t.symm}
Let $\C$ be a nonempty subset of the symmetric group $S_n$, and
let $\sigma\in \C$. Then $\C$ is a minimal AY cell if and only if
there exists a standard skew Young tableau $Q$ such that
$$
\sigma^{-1}\C=\{\pi\in S_n\,|\,Q^{\pi\inv} \hbox{ is standard}\},
$$
where $Q^{\pi\inv}$ is the tableau obtained from $Q$ by replacing each
entry $i$ by $\pi\inv(i)$.
\end{thm}

\noindent{\bf Proof.}
Given $Q$, define $f:=cont(Q)$ and use Theorem~\ref{t.cell} to conclude that
$$
\{\pi\in S_n\,|\,Q^{\pi\inv} \hbox{ is standard}\} = \C^f_{id}
$$
is a minimal AY cell containing the identity element. Thus, if $\sigma\in\C$
and $\sigma\inv\C = \C^f_{id}$ then $\C$ is a minimal AY cell.

In the other direction, if $\C$ is a nonempty minimal AY cell and $\sigma\in\C$
then $\sigma\inv\C = \C^f_{id}$ for some $\C$-generic vector $f\in V^*$.
If $L$ is the basic flat corresponding to $\C$ (see Definition~\ref{d.flat} above)
then actually $f\in L$.
Now observe that, due to the special form of the roots of $S_n$,
any (nonempty) basic flat contains a vector with integral coordinates.
By Theorem~\ref{t.5} we may thus assume that $f\in\bbz^n$, and thus
Theorem~\ref{t.tabl2} gives us the $Q$ we are looking for.
\qed

\subsection{Young Orthogonal Form}\label{s.YOF}

This subsection contains explicit representation matrices,
which are deduced from the previous analysis.
In particular, it is shown that all irreducible $S_n$-representations
may be obtained from our construction (Theorem \ref{t.irr} below).

Let ${\C}\subseteq S_n$ be a convex set and let
$f$ be an integer $\C$-generic vector.
Consider the $S_n$-representation $\T^f_{id}$.
By Corollary \ref{t.card},
a basis $\cal{B}$ of the representation space of $\T^f_{id}$
may be indexed by the set of standard Young tableaux of a certain shape.

\begin{cor}\label{t.yof2}{\bf (Young Orthogonal Form for Skew Shapes)}\\
Let $Q_0$ be a standard Young tableau of skew shape $\la/\mu$ 
($\mu$ possibly empty), 
let $f\in \bbz^n$ satisfy $\der f = \der cont(Q_0)$,
and let $\T:=\T^f_{id}$.
Then the $\T$-action of the generators of $S_n$ on the basis $\cal{B}$
is given by
$$
\T_{s_i}(v_Q)=
{1\over\hk(i)}v_Q + \sqrt{1-{1\over\hk(i)^2}}v_{Q^{s_i}}
\qquad(1\le i < n,\, v_Q\in\cal{B}),
$$
where $\hk(i):=c(i+1)-c(i)$ in $Q$, and $Q^{s_i}$ is the tableau obtained
from $Q$ by interchanging $i$ and $i+1$.
\end{cor}

\noindent{\bf Proof.} Combine Theorem~\ref{t.main12}
with Theorem~\ref{t.cell} and equation~(\ref{e.78}), noticing that,
by definition, $h(i)=\der cont(Q)_i$. \qed

\begin{thm}\label{t.irr}
Let $Q_0$ be a standard Young tableau of skew shape $\la/\mu$ 
($\mu$ possibly empty).
If $f\in \bbz^n$ satisfies $\der f = \der cont(Q_0)$ then
$$
\T^f_{id}\cong S^{\la/\mu},
$$
where $S^{\la/\mu}$ is the skew Specht module corresponding to $\la/\mu$.
\end{thm}

\noindent{\bf Proof.} For a skew shape $\la/\mu$ the
representation matrices of the generators in
Corollary~\ref{t.yof2} are identical to those given by the
classical Young orthogonal form (Theorem~\ref{t.yof-skew}). \qed

\begin{cor}\label{t.irr2}
The complete list of minimal AY representations of the symmetric
group $S_n$ is given (up to isomorphism) by the skew Specht
modules $S^{\lambda/\mu}$, where $\la/\mu$ is of order $n$
(and $\mu$ possibly empty).

In particular, every irreducible representation of the symmetric
group $S_n$ may be realized as a minimal abstract Young
representation.
\end{cor}







\section{The Irreducibles Representations of $B_n$ are AY}

We begin with the following lemma.

\begin{lem}\label{ind}
Let $(W,S)$ be a finite Coxeter system and let $J_1, J_2$ be
disjoint subsets of $S$. Let $({\cal K}_1, \rho_1)$ and $({\cal K}_2,\rho_2)$
be minimal AY pairs for $W_{J_1}$ and $W_{J_2}$, respectively,
and let $W^{J_1\cup J_2}$ be the set of representatives of minimal
length of the right cosets of the parabolic subgroup $W_{J_1\cup J_2}$ in $W$.
Then
$$
({\cal K}_1{\cal K}_2 W^{J_1\cup J_2},
(\rho_1\otimes\rho_2)\uparrow _{W_{J_1}\times W_{J_2}}^W)
$$
is a minimal AY pair for $W$.
\end{lem}

\noindent{\bf Proof.} By the definition of a minimal AY pair,
$({\cal K}_1 {\cal K}_2, \rho_1\otimes \rho_2)$ is a minimal AY
pair for $W_{J_1} W_{J_2}= W_{J_1\cup J_2}$. The lemma now follows
from Theorem~\ref{induction}.

\qed

Let $\lambda$ be a partition of $k$ ($1\le k\le n-1$),
$\mu$ a partition of $n-k$,
$P$ a standard Young tableau of shape $\lambda$ on the letters $1,\dots,k$,
and $Q$ a standard Young tableau of shape $\mu$ on the letters $k+1,\dots,n$.
Denote by $S_{[i,j]}$
the symmetric group on the letters $i, i+1, \dots, j$.

\begin{df}
A {\em shuffle} of a permutation $\pi\in S_{[1,k]}$ with
a permutation $\sigma\in S_{[k+1,n]}$ is a permutation $\tau\in S_n$
such that the letters $1,\ldots,k$ appear in $(\tau(1),\ldots,\tau(n))$
in the order $(\pi(1),\ldots,\pi(k))$ and the letters $k+1,\ldots,n$
appear in $(\tau(1),\ldots,\tau(n))$ in the order
$(\sigma(k+1),\ldots,\sigma(n))$.
\end{df}

\begin{cor}\label{t.ind-symm}
The set of all shuffles of permutations from
$$\{\pi\in S_{[1,k]}|\ P^{\pi\inv} \hbox{ is standard}\}$$
with permutations from
$$\{\sigma\in S_{[k+1,n]}|\ Q^{\sigma\inv} \hbox{ is standard}\}$$
is a minimal AY cell in $S_n$, which carries a minimal AY
representation isomorphic to the outer product
$(S^\lambda\otimes S^\mu)\uparrow_{S_k\times S_{n-k}}^{S_n}$.
\end{cor}

\noindent{\bf Proof.} Denote by $\Omega_{k,n}$ the set of all
shuffles of $(1,\ldots,k)$ with $(k+1,\ldots,n)$:
$\Omega_{k,n} = \{\tau\in S_n\,|\,\tau^{-1}(j) < \tau^{-1}(j+1), \forall j\ne k\}$.
It is well known
that $\Omega_{k,n}$ is the set of all representatives of minimal
length of right cosets of $S_{[1,k]}\times S_{[k+1,n]}$ in $S_n$.
The set considered in the corollary is the product
$$
\{\pi\in S_{[1,k]}\,|\,P^{\pi\inv} \hbox{ is standard}\}\cdot
\{\sigma\in S_{[k+1,n]}\,|\,Q^{\sigma\inv} \hbox{ is standard}\}\cdot
\Omega_{k,n}.
$$
The corollary now follows from Theorem~\ref{t.cell} and 
Theorem~\ref{t.irr} together with Lemma~\ref{ind}.

\qed

Denote the set of all shuffles considered in Corollary~\ref{t.ind-symm}
by ${\cal B}_{P,Q}$ and the associated AY representation by $\rho$.
The representation matrices of the Coxeter generators of $S_n$ on
$V_{{\cal B}_{P,Q}}$ are given by

\begin{cor}\label{t.yof-ind}
For all $1\le i< n$ and $\pi\in {\cal B}_{P,Q}$,
$$
\rho_{s_i} C_\pi=\cases{%
C_{\pi s_i}, &if either $\pi(i)\le k< \pi(i+1)$ \cr
             &or $\pi(i+1)\le k< \pi(i)$;\cr
\frac{1}{\hk(\pi,i)}C_\pi +
\sqrt{1-\frac{1}{\hk(\pi,i)^2}}C_{\pi s_i}, &otherwise,}
$$
where $\hk(\pi,i):=c(\pi(i+1))-c(\pi(i))$ in $P$ if both $\pi(i)$
and $\pi(i+1)$ are $\le k$, and in $Q$ if both are $>k$.
\end{cor}

\noindent{\bf Proof.}
Use Remark~\ref{induction1} with $W=S_n$ and
$J=\{s_1,\ldots,s_{n-1}\} \setminus \{s_k\}$.
Note that if $\pi=mr$, where $m\in S_{[1,k]}\times S_{[k+1,n]}$
and $r$ is a shuffle in $W^J = \Omega_{k,n}$, then
$rs_i\in \Omega_{k,n}$ if and only if
either $\pi^{-1}(i)\le k< \pi^{-1}(i+1)$
or $\pi^{-1}(i+1)\le k< \pi^{-1}(i)$.
The coefficients in the case $rs_i\not\in \Omega_{k,n}$
are determined, by Remark~\ref{induction1} together with
Theorem~\ref{t.tabl2} and Corollary~\ref{t.yof2}, by the content
vectors of $P$ and $Q$.

\qed

\medskip

We will show now that the subset ${\cal B}_{P,Q}\subseteq S_n$
is a minimal AY cell in $B_n$, when $S_n$ is naturally embedded in $B_n$ as
a maximal parabolic subgroup.

\begin{pro}\label{B-AY-yof}
The subset ${\cal B}_{P,Q}\subseteq B_n$ is a minimal AY
representation of $B_n$, where the action of the simple
reflections $s_i$, $1\le i<n$, is determined as in
Corollary~\ref{t.yof-ind} and the action of $s_0$ is determined by
$$
\rho_{s_0} C_\pi=\cases{%
C_{\pi},  &if $\pi(1)\le k$;\cr
-C_{\pi}, &otherwise.}
$$
\end{pro}

\noindent{\bf Proof.}
By Corollaries~\ref{t.ind-symm} and~\ref{t.yof-ind} it suffices to verify
the relations involving $s_0$. Clearly, $\rho_{s_0}$ commutes with $\rho_{s_i}$
for all $i>1$ since $\pi s_i(1) = \pi(1)$ for all $i>1$. To verify the
relation $(\rho_{s_0}\rho_{s_1})^4=1$ we have to check four cases:
\begin{itemize}
\item[(1)]
If $\pi(1)$ and $\pi(2)$ are $\le k$ then both $C_\pi$
and $C_{\pi s_1}$ are invariant under $\rho_{s_0}$;
thus $(\rho_{s_0}\rho_{s_1})^2 C_\pi = \rho_{s_1}^2 C_\pi = C_\pi$.
\item[(2)]
If $\pi(1)$ and $\pi(2)$ are $> k$ then both $C_\pi$ and $C_{\pi s_1}$
are eigenvectors of $\rho_{s_0}$ with eigenvalue $-1$;
thus again $(\rho_{s_0}\rho_{s_1})^2 C_\pi = (-\rho_{s_1})^2 C_\pi = C_\pi$.
\item[(3)]
If $\pi(1)\le k$ and $\pi(2)>k$ then
$\rho_{s_0}C_\pi = C_\pi$, $\rho_{s_0}C_{\pi s_1} = -C_{\pi s_1}$,
$\rho_{s_1}C_\pi = C_{\pi s_1}$, and $\rho_{s_1}C_{\pi s_1} = C_{\pi}$;
thus
$$
(\rho_{s_0}\rho_{s_1})^2 C_\pi= \rho_{s_0} \rho_{s_1}\rho_{s_0}
C_{\pi s_1} = -\rho_{s_0} \rho_{s_1} C_{\pi s_1}=-\rho_{s_0}
C_{\pi}=-C_{\pi}.
$$
Hence, $(\rho_{s_0}\rho_{s_1})^4 C_\pi= C_\pi$. \item[(4)] The
case $\pi(1)> k$ and $\pi(2)\le k$ is similar to Case (3) and is
left to the reader.
\end{itemize}
\qed

We deduce

\begin{thm}\label{B.irr}
All the irreducible representations of the classical Weyl group $B_n$
are minimal AY.
\end{thm}

\noindent{\bf Proof.} As before, let $\lambda$ be a partition of
$k$, $\mu$ a partition of $n-k$, $P$ a standard
Young tableau of shape $\lambda$ on the letters $1,\dots,k$  and
$Q$ a standard Young tableau of shape $\mu$ on the letters
$k+1,\dots,n$. There is a natural bijection between all pairs of
standard Young tableaux of shapes $\lambda$ and $\mu$ and elements
in the subset ${\cal B}_{P,Q}$: $(P,Q)^\pi \longleftrightarrow \pi\inv$,
where $(P,Q)^\pi$ is the ordered pair of tableaux obtained from
$(P,Q)$ by replacing each entry $i$ by $\pi(i)$. The Young Orthogonal Form
presented in Proposition~\ref{B-AY-yof} reduces to the
classical one given in Theorem~\ref{B-yof}, along this bijection.

\qed



\section{Top Elements}\label{s.inter}

This section is motivated by the following reformulation of a theorem of
Kriloff and Ram, based on results of Loszoncy~\cite{Lo}.

\begin{thm}\label{t.x2}{\rm\cite[Theorem 5.2]{KR}}
Let $W$ be a crystallographic reflection group. Then every minimal
AY cell is a left translate of an interval (in the right weak
Bruhat poset).
\end{thm}

\begin{rem}
By~\cite[Observation 3.3]{U-I}, one can assume that
the intervals are of the form $[id,w]$.
\end{rem}





%


%



\begin{df}\label{d.top}
An element $w\in W$ is a {\em top} element if the interval $[id,w]$
is a minimal AY cell which carries an {\em irreducible} AY representation
of $W$.
\end{df}

The goal of this section is to characterize the top elements in
the symmetric group $S_n$.

\begin{df}\label{d.row}
Let $Q$ be a standard skew Young tableau.
\begin{itemize}
\item[1.]
$Q$ is a {\em row (column) tableau} if and only if the
entries in each row (column) are larger than the entries in all
preceding rows (columns).
\item[2.]
The {\em row word} of $\, Q$ is obtained by reading $Q$
row by row from right to left.
The {\em column word} of $\, Q$ is obtained by reading $Q$
column by column from bottom to top.
\end{itemize}
\end{df}

\begin{exa}
The tableau
\[
\begin{array}{ccc}
1 & 2 & 3 \\
4 & 5     \\
\end{array}
\]
is a row tableau. Its row word is the permutation $[3,2,1,5,4]\in S_5$
and its column word is  $[4,1,5,2,3]\in S_5$.
\end{exa}

Our result is that there is a bijection between top elements of $S_n$ and partitions
of $n$. More explicitly:

\begin{thm}\label{t.s5}
A permutation $\pi\in S_n$ is a top element if and only if $\pi$
is the column word of a row standard Young tableau of shape $\lambda$,
where $\lambda$ is a partition of $n$.
\end{thm}

To prove this theorem we need the following lemma.

%
%
%


\begin{lem}\label{t.s7}
Let $Q$ be a standard Young tableau of order $n$.\\
The set $\{\pi\in S_n\,|\,Q^{\pi\inv} \hbox{is standard}\}$ is
an interval in the right weak Bruhat order if and only if
$Q$ is either a row tableau or a column tableau.
The maximal element in the interval is the column (respectively, row)
word of the tableau.
\end{lem}

\noindent{\bf Proof.}
Denote $B_Q:=\{\pi\in S_n|\ Q^{\pi\inv} \hbox{ is standard}\}$.

First, we prove the easy direction. Assume that $Q$ is a row
tableau. By Theorem~\ref{t.symm} and 
Proposition~\ref{t.main.convex}, $B_Q$ is convex. Thus, to prove
that $B_Q$ is an interval $[id,\pi]$ it suffices to show that
there  is a unique element $\pi\in B_Q$ such that $\pi s\not\in
B_Q$ for all $s\not\in \D(\pi)$. Indeed, if $\sigma\in B_Q$ and
$Q^{\sigma\inv}$ is not the column tableau then there exists $i$
$(1\le i\le n-1)$ such that $i+1$ is southwest of $i$ in
$Q^{\sigma\inv}$. In this case, $\sigma s_i\in B_Q$ and
$\ell(\sigma s_i)> \ell(\sigma)$. If $\sigma\in B_Q$ and
$Q^{\sigma\inv}$ is the column tableau then there is no $i$ $(1\le
i\le n-1)$ such that $i+1$ is southwest of $i$ in
$Q^{\sigma\inv}$. In this case, if $\sigma s_i\in B_Q$ then
$\ell(\sigma s_i)<\ell(\sigma)$. We conclude that the unique
maximum $\sigma$ is the column word of the row tableau $Q$.

Similarly for a column tableau $Q$.

\medskip

Now we prove the opposite direction. Assume that $Q$ is a not a
row or column tableau. We will show that $B_Q$ has at least two
maximal elements (with respect to right weak Bruhat order).
Thus, $B_Q$ is not an interval.

If the standard Young tableau  $Q$ is a not a row or column
tableau then $Q$ has at least two rows and two columns. Without
loss of generality, the letter $2$ is in the first row
(i.e., in box $(1,2)$) of $Q$.
Then the letter in box $(\la'_1,1)$ is bigger than 2.

To find the first maximal element, start with $\pi\in B_Q$ such
that $Q^{\pi\inv}$ is a column tableau
and proceed ``up" in $B_Q$.
Observe that each step in this process is a right multiplication
by $s_i$ such that the resulting permutation is in $B_Q$ and
longer. To satisfy this, $i$ and $i+1$ could not be in the same row or
column in $Q^{\pi\inv}$ (where at each step we substitute $\pi
s_i:=\pi$). Also, the letter of $Q$ in the box of $i$ in $Q^{\pi\inv}$
must be bigger than the letter of $Q$ in the box of $i+1$ in
$Q^{\pi\inv}$. The letter in $(1,2)$  is the minimal one in the
subtableau consisting all columns except the first one. Thus
cannot move along the process. We conclude

\smallskip

\noindent{\bf Claim 1.} For every $j$, $1\le j \le\la'_1+1$, the
position of the letter $j$ is invariant under this process. For
$j\le \la'_1$ the position is $(j,1)$; the position of $\la'_1+1$
is $(1,2)$.

\smallskip

Thus, we obtained one maximal element, which is determined
by a tableau with $1,...,\la'_1$ in the first column.

\smallskip

As $Q$ is not a row tableau there exists a minimal row $j$ for
which the letter in box $(j,\la_j)$ of $Q$ is bigger than the letter
in box $(j+1,1)$. To find the second maximal element, start with a
permutation in $B_Q$ determined by the standard Young tableau of the
same shape as $Q$, in which the first letters are placed in the
subshape $\{(a,b)|\ b\le \la_j\}$ in lexicographic order, and
the rest are placed in the remaining ``upper right corner" in
lexicographic order. We proceed ``up" as before; namely, by
right multiplication by $s_i$ such that the resulting permutation
is in $B_Q$ and longer.

\smallskip

\noindent{\bf Claim 2.} The set of boxes in which the letters
$1,...,j\cdot \la_j-1$  are located (in the resulting tableau) and
the locations of $j\cdot \la_j$ and $j\cdot \la_j+1$ are invariant
under this process.

\smallskip

To verify this, it suffices to show that the location of $j\cdot
\la_j+1$ is invariant under this process. Indeed, notice that as
long as $j\cdot \la_j+1$ is in the box $(j+1,1)$ $j\cdot \la_j$
must be located in box $(j,\la_j)$ (as it cannot switch with
$j\cdot\la_j-1$ which is located either in same row or column of $j\cdot\la_j$).
But the letters in boxes $(j+1,1)$ and $(j,\la_j)$ of $Q$ are
in reverse order; thus replacing $j\cdot \la_j+1$ with $j\cdot\la_j$
will shorten the permutation (and is not ``up"!).

We also cannot replace $j\cdot \la_j+1$ with $j\cdot \la_j+2$. To
verify this notice that $j\cdot \la_j+2$ has three possible
locations during our process: $(j+1,2)$, $(j+2,1)$ and $(1,\la_j+1)$.
If it is located in box $(j+1,2)$ then it is in the same row as
$j\cdot \la_j+1$, and thus switching them gives a nonstandard
tableau; thus sends us out of $B_Q$. If it is located in box $(j+2,1)$
then it is in the same column as $j\cdot \la_j+1$, and thus switching them
sends us out of $B_Q$. So, we may assume that $j\cdot \la_j+2$ is
located in box $(1,\la_j+1)$. In this case $j>1$ and, by the definition
of $j$, the letter in box $(1,\la_j+1)$ of $Q$ is less than the letter
in box $(2,1)$ of $Q$ which is, in turn,
less than the letter in box $(j+1,1)$. Thus we cannot switch $j\cdot \la_j+1$ and
$j\cdot \la_j+2$ in this case as well.

\medskip

To complete the proof, notice that the letter in box $(\la'_1,1)$ in
the first maximal tableau is $\la'_1$, while in the second maximal
tableau it is bigger. Thus the processes determine two different
maximal tableaux.
\qed

\medskip

\noindent{\bf Proof of Theorem~\ref{t.s5}.} By
Corollary~\ref{t.irr2}, the derived content vector of a standard
Young tableau of shape $\lambda/\mu$ gives an irreducible
representation if and only if $\mu$ is the empty partition. This
fact together with Theorem~\ref{t.symm} imply that $\sigma\in S_n$
is a top element if and only if
$[id,\sigma] = \{\pi\in S_n\,|\,Q^{\pi\inv} \hbox{ is standard}\}$
for some standard tableau $Q$.
Lemma~\ref{t.s7} completes the proof. \qed

\bigskip

\noindent{\bf Acknowledgments.}
The authors thank Eli Bagno and Yona Cherniavsky
for useful comments.

\end{document}